\documentclass[11pt,a4paper]{article}
\usepackage[utf8]{inputenc}
\usepackage{bm}
\usepackage[tbtags]{amsmath}
\numberwithin{equation}{section}
\usepackage{amsthm}
\usepackage{mathtools}
\usepackage{multicol}
\usepackage{relsize}
\usepackage{geometry}
\usepackage{hyperref}
\hypersetup{
    colorlinks=true,
    linkcolor=blue,
    filecolor=blue,      
    urlcolor=blue,
}
\usepackage[symbol]{footmisc} 
\title{Asymptotically optimal test for dependent multiple testing set up}
\author{Rahul Roy \footnote{Indian Statistical Institute, Mailing address: rahul\_r@isical.ac.in} , Subir Kumar Bhandari \footnote{Indian Statistical Institute, Mailing address: subir@isical.ac.in} }
\date{}

\begin{document}

\maketitle
\newtheorem{result}{Result}[section]
\newtheorem{cor}{Corollary}[result]
\makeatletter
\newcommand{\vast}{\bBigg@{4}}
\newcommand{\Vast}{\bBigg@{5}}
\makeatother
\begin{abstract}
    In this paper we explore the behaviour of dependent test statistics for testing of multiple hypothesis . To keep simplicity, we have considered a mixture normal model with equicorrelated correlation set up. With a simple linear transformation, the test statistics were decomposed into independent components, which, when conditioned appropriately generated independent variables. These were used to construct conditional tests , which were shown to be asymptotically optimal with power as large as that obtained using N.P.Lemma. We have pursued extensive simulation to support the claim.
    
\end{abstract}
\begin{center}
\end{center}
\section{Introduction:}
Multiple Hypotheses testing works as tool for maintaining the quality of overall inference rather than improving individual inferences. This is useful when drawing inference on datasets having large number of parameters. Particularly in subjects like Microarray data related to gene, astronomy, economics etc , multiple testing procedures have varied usage. Emerging almost two decades back, this relatively new direction of testing has gone through multitude of developments in quantifying the overall error for multiple testing, each with a higher degree of applicability. Primarily the idea was to restrict the type I error (false positives) under a certain level. The foremost of them was a classical method, which used the Bonferronii correction to enforce a bound on a measure called Family Wise Error Rate (FWER). Some similar methods were developed by Holm $\&$ Sture (1979), Holland et. al. (1987) ,Simes $\&$ R John(1986), Hochberg $\&$ Yosef(1988). Whereas, van der Laan, Mark J $\&$ Dudoit (2004),  Lehmann $\&$ Romano (2005) proposed multiple testing methods controlling generelised-FWER. A step up in the field of multiple testing was achieved along with the breakthrough concept of False Discovery Rate (F.D.R.) as a measure of type I error in multiple testing in their paper of Benjamini $\&$ Hochberg (1995). Sarkar (2004) provides a elaborate literature on F.D.R.
While all of those works focus mainly on independent test statistics for multiple testing, study involving dependent set up is rare. Sun $\&$ Kai (2009) uses Hidden Markov Model(H.M.M) to model and analyse a multiple testing problem.Efron(2009) addressed this through the usage of dependent Z scores.  Bogdan et al.(2011) follows a somewhat bayesian pathway and define a class of fixed threshold multiple testing rule as ABOS (assymptotically bayes optimal under sparsity) by minimizing a bayesian equivalent of FDR i.e. BFDR . Conditions for Benjamini-Hochberg procedure and Bonferroni procedures to be ABOS were established and  they approximated the corresponding random thresholds by nonrandom threshold.\\
In this paper we modify Bogdan et al.(2011) into a dependent (Multivariate Normal) setup. The conditional distribution of observations $\mathbf{X|\mu}$ and the distribution of the priors $\mathbf{\mu}$ are assumed to be equicorrelated with contrast to the independent structure in Bogdan et al.(2011). With the help of Result 1 a linear transformation breaks up the observations into linearly independent components which results in independent test statistics through conditioning. Using these independent test statistics we device a fixed threshold conditional test for the multiple testing problem. Further, following Result 2 each of these single conditional tests are shown to be asymptotically optimal (as $\tau\uparrow \infty$) by comparing the expected type II error with the corresponding conditional test obtained by N.P.Lemma. We then perform extensive simulation to validate our findings.
\section{Model specification \&  description of the problem:} \label{section2}

\subsection{Model specification:}
In  Bogdan et al.(2011) we encounter independent normal set up where both the observations and prior follow independent normal distributions. We extend the model to correlated normal both in observation and prior. We keep our model simple by choosing a simple correlation matrix (i.e. equal correlation for all the pair of variables) and equal mean (for all the variables ) . The following assumptions sum up the model:
\begin{enumerate}

    \item Conditioned on mean $\boldsymbol{\mu}$, joint distribution of $\mathbf{X}=(X_1,X_2,\cdots,X_n)$ is multivariate normal with mean $\boldsymbol{\mu}$ and covariance matrix $\sigma_\epsilon^2\Sigma_1$, where $\sigma_\epsilon^2>0$ and $\Sigma_1$ is equicorrealted correlation matrix with equal correlation $\rho_1$. i.e. $\Sigma_1=\rho_1\times\mathbf{11'}+(1-\rho_1)\times\mathrm{I}^{n\times n}$ ; $-1\leq\rho_1\leq1$ . 
    
    \item Let, $\eta_1,\eta_2,\cdots,\eta_n$ follow i.i.d \textit{Bernoulli}(p) random variables. Distribution of $\boldsymbol{\mu}$=$(\mu_1,\mu_2,\cdots,\mu_n)$ depends on the value of the unobservable (dummy) random variables  $\boldsymbol{\eta}=(\eta_1,\eta_2,\cdots,\eta_n)$ in the following way:
  
  \[ \mu_i\sim \begin{cases} N(0,\sigma_0^2) & \mathrm{if} \ \eta_i=0 \\
    N(0,\sigma_0^2+\tau^2) & \mathrm{if} \ \eta_i=1
   \end{cases} \]
   
   Here, $\sigma_0 \geq 0 $; $ \tau \uparrow \infty$ ; $ p \downarrow 0$
   \item Lastly, another crucial assumption is that, given $\boldsymbol{\eta}[=(\eta_1,\eta_2,\cdots,\eta_n)]=\boldsymbol{\eta_0}[=(\eta_{01},\eta_{02}, \cdots,\eta_{0n})] $; $\boldsymbol{\mu} $ has a multivariate normal distribution with an equicorrelated  correlation matrix $\Sigma_2$, where, $\rho_2$ is the equal correlation coefficient. i.e. $\Sigma_2=\rho_2\times\mathbf{11'}+(1-\rho_2)\times\mathrm{I}^{n\times n}$ ; $-1\leq\rho_2\leq1$. (mean and variance are obtained from the previous assumption.)
\end{enumerate}
combining these assumptions, we arrive at the final model as follows

\begin{align*}
(\eta_1,\eta_2,\cdots,\eta_n)&\stackrel{\mathrm{i.i.d}}\sim \text{Bernoulli}(p)\\
\mathrm{Given,}\boldsymbol{\eta}=\boldsymbol{\eta_{0}}, \boldsymbol{\mu}&\sim\text{MVN}(\mathbf{0},D\Sigma_2 D)
\end{align*}

where ,\[ D=\begin{bmatrix}
\sigma_1 & 0 & \cdots & 0\\
0 & \sigma_2  & \cdots & 0\\
\vdots & \vdots & \ddots& \vdots\\
0 & 0 & \cdots & \sigma_n
\end{bmatrix} \] ,
   \[ \text{with,}\sigma_i= \begin{cases}
      \sqrt{\sigma^2_0+\tau^2}   & \mathrm{if}\ \eta_{0i}=1 \\
       \sigma_0   & \mathrm{if} \ \eta_{0i}=0
    \end{cases}
 ,i=1(1)n\]
 
Therefore, given $ \boldsymbol{\eta} \ \& \ \boldsymbol{\mu},$
\begin{equation}\label{eqn1}
    \boldsymbol X\sim \text{MVN}(\boldsymbol{\mu},\sigma_\epsilon^2\Sigma_1)
\end{equation}
Marginal distribution of \textbf{\textit{X}} given $\boldsymbol{\eta}=\boldsymbol{\eta}_0$ is(i.e. we only integrate w.r.t. $\boldsymbol{\mu}$):
 \begin{equation}\label{eqn2}
\mathbf{X}\sim \text{MVN}(\mathbf{0},\sigma_\epsilon^2\Sigma_1+D\Sigma_2 D) \end{equation}
Where D is defined as before. Marginal distribution of \textbf{\textit{X}} is:

\begin{equation} \label{eqn3}
    \mathbf{X}\sim \sum_{x=0}^n\sum_{P_x}p^x(1-p)^{n-x}\text{MVN}(\mathbf{0},\sigma_\epsilon^2\Sigma_1+D_{P_x}\Sigma_2 D_{P_x})
\end{equation} 
Where, $P_x$ corresponds to a particular permutation of x 1's and (n-x) 0's. The matrix $ D_{P_x} $ is constructed as mentioned before just by replacing $\boldsymbol{\eta} $ by $P_x$.
\subsection{Description of the problem:}Here, we observe \textbf{\textit{X}} which is dependent on $\boldsymbol{\mu}$, which in turn is dependent on $\boldsymbol{\eta}$. This $\boldsymbol{\eta}$ is unobservable (dummy). Unlike in   Bogdan et al.(2011), here both of the vectors \textbf{\textit{X}} and $\boldsymbol{\mu}$ are not independent. The query naturally comes in our mind is \begin{quote}
  ``Suppose, $\boldsymbol{\eta}$ =$\boldsymbol{\eta_0}$, unobserved. We observe our data \textbf{\textit{X}} . Can we somehow get back $\boldsymbol{\eta_0}$ from this data?"
\end{quote} This simple query can be mathematically expressed in three ways.

\begin{enumerate}
    \item{\textbf{Multiple Hypothesis Testing:}}
    Let$\begin{cases}
      H_{0i}   &  \eta_{i0}=0 \\
      H_{1i}   &  \eta_{i0}=1
    \end{cases}$ , then, our problem comes down to a multiple hypothesis problem. We make a rule to reject a set of hypotheses.
    \item{\textbf{Classification problem:}} We might think this problem as a classification one, where we want to classify the data into two groups corresponding to
    \[
         \sigma_i^2 = \begin{cases}
              \sigma_0^2  \\
               \sigma_0^2 +\tau^2 
         \end{cases}  
    \].
    \item{\textbf{Subset selection Problem:}} A subset selection problem might be formulated with respect to $\sigma_i=\sqrt{\sigma_0^2+\tau^2}$.

\end{enumerate}

In this paper we shall follow the first approach. A suitable rule for accepting / rejecting the n hypotheses will be provided based on data \textbf{\textit{X}} that will eventually answer the question concerning the unobserved values $\eta_i,i=1(1)n$. \\
  We suggest the following single step fixed cut-off rule for multiple testing:
  \begin{quote}
      ``For i's with $|X_i|>K$ reject $H_{0i}$  and accept all the other $H_{0i}$'s" 
  \end{quote}
Result \ref{result2} suggests that i'th conditional test of level $\alpha$ yields K=$\mu_0+t\sqrt{\Phi_0}$. t being the solution of $\alpha=2-\Phi(t)-\Phi(t+\frac{2\mu_0}{\sqrt{\Phi_0}})$. For our problem we assume that, $\Phi_0=\sigma_\epsilon^2(1-\rho_1)+\sigma_0^2(1-\rho_2)$ is known. It can be shown that, as $p\downarrow 0$, the trimmed mean of order $\beta$,i.e., $\bar{X}_\beta$ is consistent and unbiased estimate for $\mu_0=\sigma_\epsilon Q_1+\sigma_0 Q_2$ ( $Q_1$ and $Q_2$ are among the independent components of $X_i$s as discussed in the next section). A discussion on obtaining such a rejection rule shall be included in Section \ref{section4}. A comparative study of power at fixed level of significance with the test obtained by N.P.lemma shall reveal that (see Result \ref{result2}) this test described here is asymptotically optimal as $\tau \uparrow \infty$.

\section{Independence through linear transformation:} 
The special correlation structure assumed in the observation and the prior enables us to correspond these test statistics with linear combination of independent normal variables. The following result describes the fact:

\begin{result} \label{result1}

Let $\mathbf{U}=(U_1,U_2,\cdots,U_n) \sim \text{MVN} \ (\mathbf{0},(1-\rho)\text{I}+\rho \mathbf{11'}).$ Then,

\begin{equation*}
\exists \quad \begin{rcases}
    Z_i&\stackrel{\text{i.i.d.}}{\sim}\text{N}(0,(1-\rho)) \\
    V & \sim\text{N}(0,\rho)
\end{rcases} independent   
\end{equation*}
Such that, $(U_1,U_2,\cdots,U_n)\stackrel{d}{=} (Z_1+V,Z_2+V,\cdots,Z_n+V) \quad \forall i=1(1)n$

\end{result}
 
The next corollary follows directly.

\begin{cor} \label{cor1}
Given $ V=v, U_i \stackrel {\text{i.i.d.}}{\sim} \text{N}(v,(1-\rho))$
\end{cor}
  Result \ref{result1} necessarily tells that, $\mathbf{U}$ is a \textit{oblique projection} of a scaled version of i.i.d. standard normal variables in n+1 dimensional space. Owing to this transformation, the testing procedure becomes slightly easier. Coming back to our original problem, using Result \ref{result1} it is easy to see that, given, $\boldsymbol{\eta}= \boldsymbol{\eta_0}, \exists$ 
  
 \begin{equation*}
 \begin{rcases}
     P_{1i} \stackrel{\text{i.i.d}}{\sim} & N(0,(1-\rho_1)) \quad \forall \ i=1(1)n \\
     P_{2i} \stackrel{\text{i.i.d}}{\sim}  & N(0,(1-\rho_2)) \quad  \forall \ i=1(1)n
     \\
     Q_1 \sim & N(0,\rho_1)\\
     Q_2 \sim & N(0,\rho_2)
 \end{rcases}  \ \text{independent}
\end{equation*}
 Such that,
\begin{equation} \label{eqn4}
    X_i=\sigma_\epsilon(P_{1i}+Q_1)+\sigma_i(P_{2i}+Q_2) \quad \forall i=1(1)n
\end{equation}
 Conditioning on $Q_1 \ \& \ Q_2$ we get independent random variables which are easy to deal with. From Corollary \ref{cor1} we get, 
 \begin{equation} \label{eqn5}
     (X_i|Q_1=q_1,Q_2=q_2)\stackrel{\text{\tiny{indep.}}}{\sim} N(\sigma_\epsilon q_1+\sigma_i q_2,\sigma_\epsilon^2(1-\rho_1)+\sigma_i^2(1-\rho_2))
 \end{equation} 
 \eqref{eqn5} suggests that, if we condition on $Q_1$ and $Q_2$ the observations $X_i$ become independent. So we shall get independent conditional tests if we device the tests based on $X_i|(Q_1,Q_2)$. These will be much more easier to handle than the marginal tests as, they involve an eqicorrelated correlation structure.
 \section{Using Neyman Pearson Lemma to attain asymptotic optimality:} \label{section4}
  Let us focus on the $i^{th}$ hypothesis:
  \begin{equation*}
     \begin{split}
      H_{0i}: & \ \eta_{i0}=0 \\
      H_{1i}: &\  \eta_{i0}=1
     \end{split} 
  \end{equation*}
Following  \eqref{eqn5} which is equivalent to:
  \begin{equation*}
  \begin{split}
       H_{0i}: & \ (X_i|Q_1=q_1,Q_2=q_2)\sim N(\sigma_\epsilon q_1+\sigma_0 q_2,\sigma_\epsilon^2(1-\rho_1)+\sigma_0^2(1-\rho_2)) \\
      H_{1i}: & \ (X_i|Q_1=q_1,Q_2=q_2)\sim N(\sigma_\epsilon q_1+\sigma_\tau q_2 ,\sigma_\epsilon^2(1-\rho_1)+\sigma_\tau^2(1-\rho_2))
\end{split}
\end{equation*}
  where, $\sigma_\tau^2=\sigma_0^2+\tau^2$. Note that, under $H_{0i}$, $E(X_i|Q_1=q_1,Q_2=q_2)=\sigma_\epsilon q_1+\sigma_0 q_2$ and $var(X_i|Q_1=q_1,Q_2=q_2)=\sigma_\epsilon^2(1-\rho_1)+\sigma_0^2(1-\rho_2)) $, both of which are constant irrespective of the value of $\tau$. However, under $H_{1i}$, $E(X_i|Q_1=q_1,Q_2=q_2)=\sigma_\epsilon q_1+\sqrt{\sigma_0^2+\tau^2}q_2$ ,which tend to $\infty$ or $-\infty$ depending on the sign of $q_2$, and $var(X_i|Q_1=q_1,Q_2=q_2)=\sigma_\epsilon^2(1-\rho_1)+(\sigma_0^2+\tau^2)(1-\rho_2)) $ which tends to $\infty$ with $\tau\uparrow \infty$. Therefore, it is expected that,under $H_{1i}$, $|X_i|$ values are higher than that of under $H_{0i}$. Hence, our form of test discussed in section \ref{section2} is justified. To justify analytically let us define, $Y_i=(\frac{X_i}{\tau}|Q_1=q_1,Q_2=q_2).$
  Then, 
  \begin{equation*}
      \begin{split}
      H_{0i}:& \ Y_i\sim N(\frac{\mu_0}{\tau},\frac{\Phi_0}{\tau^2})  \\
      H_{1i}: & \ Y_i \sim N(\mu_1(\tau),\Phi_1(\tau))
      \end{split}
  \end{equation*} 
 
  Where,
  \begin{align*}
  \mu_0=&\sigma_\epsilon q_1+\sigma_0 q_2\\
   \Phi_0=&\sigma_\epsilon^2(1-\rho_1)+\sigma_0^2(1-\rho_2)\\
  \mu_1(\tau)=&\frac{\sigma_\epsilon q_1}{\tau}+ q_2\sqrt{\big(\frac{\sigma_0}{\tau}\big)^2+1}\\ \Phi_1(\tau)=&\frac{\sigma_\epsilon^2(1-\rho_1)}{\tau^2}+(\big(\frac{\sigma_0}{\tau}\big)^2+1)(1-\rho_2)
  \end{align*}
  
  Note that: Given, $Q_1=q_1 \ \& \ Q_2=q_2 ,\ \mu_0, \Phi_0$ are finite and \[ \mu_1(\tau) \xrightarrow{\tau \uparrow \infty} q_2 \ \& \ \Phi_1(\tau) \xrightarrow{\tau \uparrow \infty} (1-\rho_2).\]
  \paragraph{\textbf{Rationale behind taking $Y_i=(\frac{X_i}{\tau}|Q_1=q_1,Q_2=q_2)$ for analysis purpose:}}
  \begin{enumerate}
  \item As mentioned earlier, we get advantage on conditioning as the observations can be expressed as linear combinations of independent variables. Thus, the test statistics chosen for testing different hypotheses become independent this way.
  
      \item As $\tau \uparrow \infty,$ under $H_{1i}, \ \text{E}(X_i)\uparrow \infty$ or $\text{E}(X_i)\downarrow -\infty$ \& \ $\text{Var}(X_i)\uparrow \infty$  but those for $Y_i$ goes to constant. So handling the variable under alternative remain under control. Also $\tau $ being a positive value, direction of the rejection region based on $Y_i$ and $(X_i|Q_1=q_1,Q_2=q_2)$ are same
      
  \end{enumerate}
  
  \textbf{N-P lemma} says,the most powerful test of its size rejects  $H_{0i}$ iff $\frac{f_1(y_i)}{f_0(y_i)}>k$ , where $f_0(.) \  \& \ f_1(.)$ are respectively the null and alternative distributions of $Y_i$.
  Since,  Y is dependent on $\tau,q_1,q_2,$ This constant k is possibly dependent on $\tau,q_1,q_2$ and is denoted as $k_1(\tau,q_1,q_2)$.
  So, N.P. lemma $\implies \ H_{0i}$ is rejected if,
  $$\frac{\frac{1}{\sqrt{2\pi\Phi_1(\tau)}}\exp\bigg[-\frac{(y_i-\mu_1(\tau))^2}{2\Phi_1(\tau)} \bigg]}{\frac{\tau}{\sqrt{2\pi\Phi_0}}\exp\bigg[ -\frac{\big (y_i-\frac{\mu_0}{\tau}\big)^2\tau^2}{2\Phi_0}\bigg]}>k_1(\tau,q_1,q_2) $$
  i.e.
  $$e^{\bigg( \frac{y_i^2\tau^2}{2\Phi_0}-\frac{y_i^2}{2\Phi_1(\tau)}-\frac{2y_i\mu_0\tau}{2\Phi_0}+\frac{2y_i\mu_1(\tau)}{2\Phi_1(\tau)}+\frac{\mu_0^2}{2\Phi_0}- \frac{\mu_1^2(\tau)}{2\Phi_1(\tau)}-\log\big[\tau\sqrt{\Phi_1(\tau)} \big]-\log\sqrt{}\Phi_0\bigg)}>k_1(\tau,q_1,q_2) $$
  i.e.
  $$\frac{y_i^2\tau^2}{2\Phi_0}+o(\tau^2)>\log(k_1(\tau,q_1,q_2)) $$
  Since, $\tau \uparrow \infty,$ we ignore rest of the terms which are of small order of $\tau^2$ and thus our testing ,method becomes: \begin{quote}
      Given $Q_1=q_1, Q_2=q_2$ reject $H_{0i}$ if $Y_i^2>k_2(\tau,q_1,q_2)$ i.e. if $|X_i|>K(\tau,q_1,q_2)$
  \end{quote}
  Here's how we get the test stated previously. Note that, we have ignored all the terms having lower powers of $\tau$. Surely, as a result, we shall get skewed conclusions compared to that obtained from the N.P.Lemma for small values of Tau. And since NP. lemma yields the most powerful test, power obtained by this test will be somewhat lesser.   
     Our job is however to find the efficiency of the method at large values of $\tau$. In the next result, we shall prove that, for large $\tau$, this difference in type II error (and hence power) is infinitesimal up to a certain order of $\tau$.  
  \begin{result} \label{result2}
  For a given level $\alpha$ (with $\alpha$ small),
  \begin{enumerate}
      \item The level $\alpha$ conditional test of the form $|X_i|>K(\tau,q_1,q_2)$ has  $K(\tau,q_1,q_2)=\mu_0+t\sqrt{\Phi_0}$. t being the solution of $\alpha=2-\Phi(t)-\Phi(t+\frac{2\mu_0}{\sqrt{\Phi_0}})$.
      \item For $\tau \uparrow \infty$ the level $\alpha$ test described by N.P.Lemma has rejection region of the form:
     
     \[
  \{\text{X}:\text{X}> K_1^{'\alpha}(\tau)\}\cup\{y:y< K_2^{'\alpha}(\tau)\}
  \]
     where 
    \begin{equation*}
       \begin{split}
           K_1^{'\alpha}(\tau)=& (\mu_0+\sqrt{\Phi_0}z_{\frac{\alpha}{2}}) +_o(1)
       \end{split}
   \end{equation*}
     And ,  
      \begin{equation*}
       \begin{split}
           K_2^{'\alpha}(\tau)=& (\mu_0-\sqrt{\Phi_0}z_{\frac{\alpha}{2}}) +_o(1)
       \end{split}
   \end{equation*}
   Where, $z_\gamma$ is the upper 100$\gamma\%$ upper quantile of standard normal distribution.
      \item The expected type II error of the conditional test followed by the N.P. lemma is of order $\frac{1}{\tau}$ and is given by $\frac{1}{\tau}\frac{2\sqrt{\Phi_0}z_{\frac{\alpha}{2}}}{\sqrt{2\pi}}$.
      \item For sufficiently large $\tau$, the difference between the expected type II errors of the conditional tests described in this paper and that corresponding to the N.P. Lemma goes to 0 faster than $\frac{1}{\tau}$.
  \end{enumerate}
  
  \end{result}
   
   \paragraph{Remarks:}
   \begin{enumerate}
       \item The 3rd part of the result is suggestive that the power of the  most powerful test , i.e. the test based on N.P. lemma , goes to 1 at a rate of $\frac{1}{\tau}$. The 4th part shows that, sufficiently large $\tau$ enables the power of the two tests to differ at most at a rate less than $\frac{1}{\tau}$. Now, since the test based on N.P. yields the most powerful test, asymptotically, our test acquire the same attributes too. Therefore our test happens to be the asymptotically optimal test for the given hypothesis testing problem (for large $\tau$). 
       \item  From the first part of the result, we get, the level $\alpha $ test for ith hypothesis is: $|X_i|>\mu_0+t\sqrt{\Phi_0}$. Where, $\mu_0=\sigma_\epsilon q_1+\sigma_0 q_2$, $\Phi_0=\sigma_\epsilon^2(1-\rho_1)+\sigma_0^2(1-\rho_2)$ and t is the solution of: $\alpha=2-\Phi\bigg(\frac{2\mu_0}{\sqrt{\Phi_0}}+t\bigg)-\Phi(t)$. Since, this cut off value is not dependent on i, we note down which of the $|X_i|$ values exceed this cut off and reject the corresponding null hypotheses. We accept all the other null hypotheses. 
       \item For now we assume that, $\sigma_\epsilon, \sigma_0,\rho_1\&\rho_2$ are all known. i.e. we know the value of $\Phi_0$.
       \item It can be shown that, under the assumption of $p\downarrow 0$,the trimmed mean of order $\beta$ ($\bar{X}_\beta$) is an unbiased and consistent estimator of $\mu_0$. Hence , we estimate $\mu_0$ by $\bar{X}_\beta$
       
   \end{enumerate}
   \section{Simulation Study:}
     From the proof of Result \ref{result2} we get the expected type II error of our method as of the form:$\frac{1}{\tau}\frac{2\sqrt{\Phi_0}z_{\frac{\alpha}{2}}}{\sqrt{2\pi}}$ with an error atmost of order $_o(\frac{1}{\tau})$  . We verify this result in this section via extensive simulation study. The following simulation scheme has been adopted.
     \paragraph{Simulation Scheme:}
     \begin{enumerate}
         \item Fix n=500. Now generate a n variate vector of observation based on the model we have mentioned in the begining. 
         \item Compute trimmed-mean $\bar{X}_\beta$ for $\beta=0.05$. Subtract this quantity from all the observations. Now we compute the cut off value as the 100$\alpha\%$ upper quantile of the absolute values of the transformed observations which follows the null distribution. Call this K
         \item Reject those null hypotheses for which, absolute values of the transformed observation exceeds K. Repeat the same for 500 times and take the average  proportion of false positive, proportion of false negative  and their standard deviations.
         \item  Compute corresponding expected value of type II error as obtained by Result \ref{result2}. Compare the proportions of false positives with this values.
     \end{enumerate}
     Here, We assumed, $\sigma_\epsilon=\sigma_0=1,$ We change $\rho_1 = \rho_2=\rho$ in \{0,0.1,0.4,0.7\} and we conduct the whole process for $\tau=\{1,3,7,15,30,50,100\}$. These were repeated for p=0.05 after p=0.1
     
     We observe that proportion of false negatives decreases with increasing $\tau$. When $\rho_1=\rho_2=\rho$ is small, the difference between the proportion of false negatives and the corresponding expected Type II errors are negligible for $\tau$ as large as 100. However , as  $\rho$ increases, there is a visible difference , which occurs mainly because we ignored the terms involving smaller power of $\tau$, which is expected to vanish if we take higher values of $\tau$.
        The simulation scheme requires that $\alpha<$p. We considered $\alpha=0.05$. However, result was satisfactory, even when $\alpha=p=0.05$.

\begin{table}[!ht]
    \centering
    \caption{n=500 , p=0.1 , $\rho_1$=$\rho_2$=0}
    \begin{tabular}{|c|ccccc|}
    \hline
      $\tau$ & p.f.p. & sd(p.f.p) & p.f.n. & sd(p.f.n)& E(Type II error)\\
    \hline  
       1   & 0.050960839 & 2.69$\times10^{-05}$ & 0.890047716 & 0.002030549 & 2.211582528 \\
3   & 0.050918769 & 2.66$\times10^{-05}$ & 0.588603949 & 0.0034967   & 0.737194176 \\
7   & 0.05098757  & 2.57$\times10^{-05}$ & 0.300295167 & 0.003114475 & 0.315940361 \\
15  & 0.050982728 & 2.60$\times10^{-05}$ & 0.144813067 & 0.002264286 & 0.147438835 \\
30  & 0.050974072 & 2.67$\times10^{-05}$ & 0.072412769 & 0.001599078 & 0.073719418 \\
50  & 0.050958139 & 2.65$\times10^{-05}$ & 0.048803314 & 0.001425052 & 0.044231651 \\
100 & 0.050988826 & 2.62$\times10^{-05}$ & 0.022888149 & 0.000991919 & 0.022115825 \\
\hline
\end{tabular} 
\end{table}

\begin{table}[!ht]
\centering
\caption{n=500 , p=0.1 , $\rho_1$=$\rho_2$=0.1}
\begin{tabular}{|c|ccccc|}
\hline
      $\tau$ & p.f.p. & sd(p.f.p) & p.f.n. & sd(p.f.n)& E(Type II error)\\
    \hline 
1   & 0.05095194  & 2.68$\times10^{-05}$ & 0.886153991 & 0.002213726 & 2.098091407 \\
3   & 0.050956782 & 2.61$\times10^{-05}$ & 0.588264374 & 0.003283132 & 0.699363802 \\
7   & 0.050903549 & 2.73$\times10^{-05}$ & 0.296977689 & 0.003065582 & 0.299727344 \\
15  & 0.050985516 & 2.83$\times10^{-05}$ & 0.143202257 & 0.00228495  & 0.13987276  \\
30  & 0.050910616 & 2.69$\times10^{-05}$ & 0.075503167 & 0.001820517 & 0.06993638  \\
50  & 0.050990047 & 2.62$\times10^{-05}$ & 0.044254001 & 0.001261119 & 0.041961828 \\
100 & 0.050938893 & 2.72$\times10^{-05}$ & 0.023520354 & 0.000993864 & 0.020980914\\
\hline
\end{tabular}
\end{table}
\begin{table}[!ht]
\centering
\caption{n=500 , p=0.1 , $\rho_1$=$\rho_2$=0.4}
\begin{tabular}{|c|ccccc|}
\hline
      $\tau$ & p.f.p. & sd(p.f.p) & p.f.n. & sd(p.f.n)& E(Type II error)\\
    \hline
1   & 0.05098732  & 2.70$\times10^{-05}$ & 0.879330195 & 0.002660675 & 1.71308446  \\
3   & 0.050945272 & 2.73$\times10^{-05}$ & 0.554261255 & 0.005605063 & 0.571028153 \\
7   & 0.050974659 & 2.71$\times10^{-05}$ & 0.269348361 & 0.00414948  & 0.244726351 \\
15  & 0.050935446 & 2.78$\times10^{-05}$ & 0.128871825 & 0.002572435 & 0.114205631 \\
30  & 0.051006447 & 2.72$\times10^{-05}$ & 0.068629977 & 0.001877376 & 0.057102815 \\
50  & 0.050939078 & 2.78$\times10^{-05}$ & 0.038717404 & 0.001333223 & 0.034261689 \\
100 & 0.050975602 & 2.70$\times10^{-05}$ & 0.020250571 & 0.00094914  & 0.017130845\\
\hline
\end{tabular}
\end{table}
\begin{table}[!ht]
\centering
\caption{n=500 , p=0.1 , $\rho_1$=$\rho_2$=0.7}
\begin{tabular}{|c|ccccc|}
\hline
      $\tau$ & p.f.p. & sd(p.f.p) & p.f.n. & sd(p.f.n)& E(Type II error)\\
    \hline
1   & 0.050992441 & 2.68$\times10^{-05}$ & 0.851410397 & 0.003469421 & 1.211333639 \\
3   & 0.050939205 & 2.62$\times10^{-05}$ & 0.472698246 & 0.009387919 & 0.40377788  \\
7   & 0.050945869 & 2.68$\times10^{-05}$ & 0.217796774 & 0.006142334 & 0.173047663 \\
15  & 0.050985557 & 2.81$\times10^{-05}$ & 0.101505957 & 0.003574137 & 0.080755576 \\
30  & 0.051001283 & 2.75$\times10^{-05}$ & 0.055259237 & 0.002201328 & 0.040377788 \\
50  & 0.051010085 & 2.65$\times10^{-05}$ & 0.033167118 & 0.001480051 & 0.024226673 \\
100 & 0.051012557 & 2.77$\times10^{-05}$ & 0.01597581  & 0.000918522 & 0.012113336\\
\hline
\end{tabular}
\end{table}
\footnote{p.f.p.=proportion of false positives, p.f.n.=proportion of false positives, s.d.=standard deviation}
\begin{table}[!ht]
\centering
\caption{n=500 , p=0.05 , $\rho_1$=$\rho_2$=0}
\begin{tabular}{|c|ccccc|}
\hline
      $\tau$ & p.f.p. & sd(p.f.p) & p.f.n. & sd(p.f.n)& E(Type II error)\\
    \hline
1   & 0.050665746 & 2.39$\times10^{-05}$ & 0.892008958 & 0.002858644 & 2.211582528 \\
3   & 0.05067558  & 2.36$\times10^{-05}$ & 0.587653248 & 0.004701227 & 0.737194176 \\
7   & 0.050687026 & 2.41$\times10^{-05}$ & 0.300371923 & 0.004113348 & 0.315940361 \\
15  & 0.050708276 & 2.37$\times10^{-05}$ & 0.143438729 & 0.003151845 & 0.147438835 \\
30  & 0.050697258 & 2.35$\times10^{-05}$ & 0.075758958 & 0.002447155 & 0.073719418 \\
50  & 0.050746188 & 2.39$\times10^{-05}$ & 0.04310378  & 0.001931703 & 0.044231651 \\
100 & 0.050709201 & 2.39$\times10^{-05}$ & 0.019661177 & 0.001247551 & 0.022115825\\
\hline
\end{tabular}
\end{table}
\newpage
\begin{table}[!ht]
\centering
\caption{n=500 , p=0.05 , $\rho_1$=$\rho_2$=0.1}
\begin{tabular}{|c|ccccc|}
\hline
      $\tau$ & p.f.p. & sd(p.f.p) & p.f.n. & sd(p.f.n)& E(Type II error)\\
    \hline
1   & 0.050728136 & 2.40$\times10^{-05}$  & 0.88493528  & 0.003003929 & 2.098091407 \\
3   & 0.050717412 & 2.47$\times10^{-05}$  & 0.588207355 & 0.004532271 & 0.699363802 \\
7   & 0.050697125 & 2.50$\times10^{-05}$  & 0.30168262  & 0.004457361 & 0.299727344 \\
15  & 0.050726196 & 2.27$\times10^{-05}$  & 0.143843566 & 0.00323125  & 0.13987276  \\
30  & 0.050702404 & 2.31$\times10^{-05}$  & 0.075416828 & 0.002399414 & 0.06993638  \\
50  & 0.050668561 & 2.34$\times10^{-05}$  & 0.041907941 & 0.001871762 & 0.041961828 \\
100 & 0.050700848 & 2.31$\times10^{-05}$  & 0.022051089 & 0.001353731 & 0.020980914\\
\hline
\end{tabular}
\end{table}
\begin{table}[!ht]
\centering
\caption{n=500 , p=0.05 , $\rho_1$=$\rho_2$=0.4}
\begin{tabular}{|c|ccccc|}
\hline
      $\tau$ & p.f.p. & sd(p.f.p) & p.f.n. & sd(p.f.n)& E(Type II error)\\
    \hline
1   & 0.050720465 & 2.45$\times10^{-05}$ & 0.875870271 & 0.003182193 & 1.71308446  \\
3   & 0.050689452 & 2.35$\times10^{-05}$ & 0.558125836 & 0.006094067 & 0.571028153 \\
7   & 0.050692737 & 2.35$\times10^{-05}$ & 0.270845222 & 0.005021932 & 0.244726351 \\
15  & 0.050680071 & 2.35$\times10^{-05}$ & 0.130391194 & 0.003396935 & 0.114205631 \\
30  & 0.050725066 & 2.31$\times10^{-05}$& 0.064686563 & 0.002320392 & 0.057102815 \\
50  & 0.050717542 & 2.37$\times10^{-05}$ & 0.039116386 & 0.001881571 & 0.034261689 \\
100 & 0.050664361 & 2.24$\times10^{-05}$ & 0.018966162 & 0.001305745 & 0.017130845\\
\hline
\end{tabular}
\end{table}
\begin{table}[!ht]
\centering
\caption{n=500 , p=0.05 , $\rho_1$=$\rho_2$=0.7}
\begin{tabular}{|c|ccccc|}
\hline
      $\tau$ & p.f.p. & sd(p.f.p) & p.f.n. & sd(p.f.n)& E(Type II error)\\
    \hline
1   & 0.050727688 & 2.34$\times10^{-05}$ & 0.856534072 & 0.003951657 & 1.211333639 \\
3   & 0.050662938 & 2.32$\times10^{-05}$ & 0.481570553 & 0.009932291 & 0.40377788  \\
7   & 0.050724644 & 2.35$\times10^{-05}$ & 0.217786875 & 0.006582882 & 0.173047663 \\
15  & 0.050754802 & 2.37$\times10^{-05}$ & 0.104966674 & 0.003909864 & 0.080755576 \\
30  & 0.050684963 & 2.33$\times10^{-05}$ & 0.053428116 & 0.002627073 & 0.040377788 \\
50  & 0.050667026 & 2.38$\times10^{-05}$ & 0.033783384 & 0.001884929 & 0.024226673 \\
100 & 0.050714643 & 2.33$\times10^{-05}$ & 0.015958144 & 0.001191515 & 0.012113336\\
\hline

\end{tabular}
\end{table}
\clearpage
\section{Proofs:}
 \subsection{Proof of  Result \ref{result1}}
     \begin{proof}
     \begin{equation}
     \begin{split}
    \text{MGF of \textbf{U}}=&(U_1,U_2,\cdots,U_n):\\
    M_{ \text{\textbf{U}}}(\text{\textbf{t}})=&\text{E}(\exp(\text{\textbf{t$^\prime$U}}))= \text{E}\bigg[\exp\bigg (\sum_{i=1}^nt_iU_i\bigg )\bigg]\\
     =&\exp\bigg (\frac{1}{2}\bigg[(1-\rho)\text{\textbf{t$^\prime$t}}+\rho\text{\textbf{t$^\prime$11$^\prime$t}}\bigg]\bigg) \\
     =&\exp\bigg(\frac{1}{2}\bigg[(1-\rho)\sum_{i=1}^nt_i^2+\rho(\sum_{i=1}^nt_i)^2 \bigg] \bigg) \\
    =&\exp\bigg(\frac{1}{2}\rho(\sum_{i=1}^nt_i)^2 \bigg)\times\prod_{i=1}^n\exp\bigg(\frac{1}{2}(1-\rho)t_i^2 \bigg)\\
     \end{split}
      \end{equation} \label{eqn6}
     $\exists \quad \begin{array}{cc}
    Z_i \stackrel{\text{i.i.d.}}{\sim}  & \text{N}(0,(1-\rho)) \\
    V \sim & \text{N}(0,\rho)
\end{array} \bigg \} \text{independent,} $ with , $  
    \text{M}_{Z_i}(t)= \text{E}(\exp(tZ_i)) = $\\ $\exp\big(\frac{1}{2}(1-\rho) t^2\big)\ \forall i =1(1)n$ $\&$
   $ \text{M}_{V}(t)= \text{E}(\exp(tV))= \exp\big(\frac{1}{2}\rho t^2\big)
  $
  Therefore, \ref{eqn6} $\implies$
  \begin{equation*}
  \begin{split}
   M_{ \text{\textbf{U}}}\mathbf{(t)}=& \text{E}( \exp\bigg[ V\bigg(\sum_{i=1}^nt_i\bigg)\bigg])\times \prod_{i=1}^n\text{E}[\exp(Z_it_i ] )\\
  =&\text{E}(\exp\bigg[\sum_{i=1}^n(t_i(V+Z_i)) \bigg])\\
  \text{So}, \ \text{E}(\exp\bigg (\sum_{i=1}^nt_iU_i\bigg ))=& \text{E}(\exp\bigg[\sum_{i=1}^n(t_i(V+Z_i)) \bigg])\\
  \text{i.e.} \ \text{M}_{\text{\textbf{U}}}(\mathbf{t})=&\text{M}_{\text{\textbf{Z}}+V\mathbf{1}}(\mathbf{t})\\
  \end{split}
\end{equation*}
  [Here, $\text{\textbf{Z}}=(Z_1,Z_2,\cdots,Z_n)^\prime \ \& \ \text{\textbf{U}}=(U_1,U_2,\cdots,U_n)^\prime  \ \& \ \text{\textbf{1}}=(1,1,\cdots,1)^\prime$]
\[
\implies (U_1,U_2,\cdots,U_n)\stackrel{d}{=} (Z_1+V,Z_2+V,\cdots,Z_n+V) \quad \forall i=1(1)n
\]

 \end{proof}
 \subsection{Proof of  Result \ref{result2}}
   \begin{proof}
   Here, we restate the hypotheses once again:
    \begin{equation*}
        \begin{split}
      H_{0i}:& Y_i\sim N(\frac{\mu_0}{\tau},\frac{\Phi_0}{\tau^2})  \\
      H_{1i}: & Y_i \sim N(\mu_1(\tau),\Phi_1(\tau))   
        \end{split}
    \end{equation*} 
   Where,
  \begin{align*}
  \mu_0=&\sigma_\epsilon q_1+\sigma_0 q_2\\
   \Phi_0=&\sigma_\epsilon^2(1-\rho_1)+\sigma_0^2(1-\rho_2)\\
  \mu_1(\tau)=&\frac{\sigma_\epsilon q_1}{\tau}+\sqrt{\big(\frac{\sigma_0}{\tau}\big)^2+1} q_2\\ \Phi_1(\tau)=&\frac{\sigma_\epsilon^2(1-\rho_1)}{\tau^2}+(\big(\frac{\sigma_0}{\tau}\big)^2+1)(1-\rho_2)
  \end{align*}
  We first fix a level of significance $\alpha$ (small).The testing procedure we suggested will reject $H_{0i}$ if, $|X_i|>K(q_1,q_2)$. Now let's find out the type II error obtained by this test at a level of significance $\alpha$. The size of the test would be:(Here, \ $P_{H_0}$ stands for probability under $H_0$ )
   \begin{equation} \label{eqn7}
   \begin{split}
  \text{Size}=& P_{H_0}(|X_i|>K(q_1,q_2)) \\ 
    = & P_{H_0}\big (|Y_i|>\frac{K(q_1,q_2)}{\tau} \big ) \\
    =&P_{H_0}\big (Y_i>\frac{K(q_1,q_2)}{\tau}\big )+P_{H_0}\big (Y_i<-\frac{K(q_1,q_2)}{\tau}\big )\\
    =& P\big (Z>\frac{\frac{K(q_1,q_2)-\mu_0}{\tau}}{\frac{\sqrt{\Phi_0}}{\tau}}\big)+P\big (Z<-\frac{\frac{K(q_1,q_2)+\mu_0}{\tau}}{\frac{\sqrt{\Phi_0}}{\tau}}\big) \\
    & (\text{Where, } Z\sim N(0,1) \text{ and independent with Q}_1 \text{ and Q}_2 \text{ under H}_{0i})\\
   =&P\big (Z>\frac{K(q_1,q_2)-\mu_0}{\sqrt{\Phi_0}}\big)+P\big (Z< -\frac{K(q_1,q_2)+\mu_0}{\sqrt{\Phi_0}}\big) \\ 
   =&1-\Phi\big(\frac{K(q_1,q_2)-\mu_0}{\sqrt{\Phi_0}}\big)+\Phi\big(-\frac{K(q_1,q_2)+\mu_0}{\sqrt{\Phi_0}} \big) \\
   =&1-\Phi(t_1)+\Phi(t_2) \\
   \end{split}
   \end{equation}
    Where,
    \begin{align*}
       t_1=&\frac{K(q_1,q_2)-\mu_0}{\sqrt{\Phi_0}}\text{ and,}\\ t_2=&-\frac{K(q_1,q_2)+\mu_0}{\sqrt{\Phi_0}}\\ 
   \implies t_1+t_2=&-\frac{2\mu_0}{\sqrt{\Phi_0}}\\
  \implies t_2=&-\frac{2\mu_0}{\sqrt{\Phi_0}}-t_1
    \end{align*} 
    Therefore, \ref{eqn7} $\implies$\begin{equation*}
        \begin{split}
            \text{Size}=&1-\Phi(t_1)+\Phi(-\frac{2\mu_0}{\sqrt{\Phi_0}}-t_1)\\ 
              =&2-\Phi(t_1)-\Phi(\frac{2\mu_0}{\sqrt{\Phi_0}}+t_1) 
        \end{split}
    \end{equation*}
   equating Size to $\alpha$, we get, 
   \[
   \alpha=2-\Phi(t_1)-\Phi(\frac{2\mu_0}{\sqrt{\Phi_0}}+t_1) 
   \]
    Now recall, $t_1=\frac{K(q_1,q_2)-\mu_0}{\sqrt{\Phi_0}} $
  \begin{equation}
  \implies K(q_1,q_2)=\sqrt{\Phi_0}t_1+\mu_0 
  \end{equation}
    So, if $t$ is the solution of: $\alpha=2-\Phi(t)-\Phi(\frac{2\mu_0}{\sqrt{\Phi_0}}+t) $ , we shall reject $H_{0i}$ at $\alpha $ level of significance if $|X_i|>\sqrt{\Phi_0}t+\mu_0 $. This proves the first  part of Result \ref{result2}
    .
   
   Therefore:
   \begin{equation} \label{eqn8}
   \begin{split}
   \text{Power}=&P_{H_{1i}}(|X_i|>\sqrt{\Phi_0}t+\mu_0 )\\
    =&P_{H_{1i}}\big(|Y_i|>\frac{\sqrt{\Phi_0}t+\mu_0}{\tau}\big)\\
   =&P_{H_{1i}}\big (Y_i>\frac{\sqrt{\Phi_0}t+\mu_0}{\tau}\big )+P_{H_{1i}}\big( Y_i<-\frac{\sqrt{\Phi_0}t+\mu_0}{\tau}\big)\\
   =&P\big (Z>\frac{\frac{\sqrt{\Phi_0}t+\mu_0}{\tau}-\mu_1(\tau)}{\sqrt{\Phi_1(\tau)}}\big )+P\big( Z<\frac{-\frac{\sqrt{\Phi_0}t+\mu_0}{\tau}-\mu_1(\tau)}{\sqrt{\Phi_1}(\tau)}\big)\\
    & (\text{Where,} Z\sim N(0,1) \text{ and independent with Q}_1 \text{ and Q}_2 \text{ under H}_{1i})\\
   =&1-\Phi\bigg(\frac{\frac{\sqrt{\Phi_0}t+\mu_0}{\tau}-\mu_1(\tau)}{\sqrt{\Phi_1(\tau)}} \bigg)+\Phi\bigg(\frac{-\frac{\sqrt{\Phi_0}t+\mu_0}{\tau}-\mu_1(\tau)}{\sqrt{\Phi_1(\tau)}} \bigg)\\
   \end{split}
   \end{equation}
   Now let's look at the following quantities:
  
 \begin{equation*}
     \begin{split}
 \mu_1(\tau)=&\frac{\sigma_\epsilon q_1}{\tau}+q_2\sqrt{\frac{\sigma_0^2}{\tau^2}+1}\\
 =&\frac{\sigma_\epsilon q_1}{\tau}+q_2(1+\frac{\sigma_0^2}{2\tau^2}+\cdots)\\
 =& q_2+\frac{\sigma_\epsilon q_1}{\tau}+_o(\frac{1}{\tau})\\
 \end{split}
 \end{equation*}
 
 \begin{equation*}
 \begin{split}
  \Phi_1(\tau)=&\frac{\sigma_\epsilon^2(1-\rho_1)}{\tau^2}+\big( \frac{\sigma_0^2}{\tau^2}+1\big)(1-\rho_2)\\
  =&(1-\rho_2)+\frac{\sigma_\epsilon^2(1-\rho_1)+\sigma_0^2(1-\rho_2)}{\tau^2}\\
  =&(1-\rho_2)+\frac{\Phi_0}{\tau^2}\\
  \end{split}
  \end{equation*}
  
  \begin{equation*}
      \begin{split}
  \implies \sqrt{\Phi_1(\tau)}=&\sqrt{(1-\rho_2)}\big(1+\frac{\Phi_0}{(1-\rho_2)\tau^2} \big)^{0.5}\\
  =&\sqrt{(1-\rho_2)}\big( 1+\frac{\Phi_0}{2(1-\rho_2)\tau^2}+\cdots\big)\\
  =&\sqrt{(1-\rho_2)}+_o(\frac{1}{\tau})\\
  \end{split}
  \end{equation*}
  
  \paragraph{}
  Hence,from \ref{eqn8} we have; Power
  \begin{equation} \label{eqn9}
 \begin{split}
 =&1-\Phi\bigg(\frac{\frac{\sqrt{\Phi_0}t+\mu_0}{\tau}-q_2-\frac{\sigma_\epsilon q_1}{\tau}-_o(\frac{1}{\tau})}{\sqrt{(1-\rho_2)}+_o(\frac{1}{\tau})} \bigg) +\Phi\bigg(\frac{-\frac{\sqrt{\Phi_0}t+\mu_0}{\tau}-q_2-\frac{\sigma_\epsilon q_1}{\tau}-_o(\frac{1}{\tau})}{\sqrt{(1-\rho_2)}+_o(\frac{1}{\tau})} \bigg) \\
  =&1-\Phi\bigg(\frac{\frac{\sqrt{\Phi_0}t+\mu_0-\sigma_\epsilon q_1}{\tau}-q_2}{\sqrt{(1-\rho_2)}}+_o(\frac{1}{\tau}) \bigg)+\Phi\bigg(-\frac{\frac{\sqrt{\Phi_0}t+\mu_0+\sigma_\epsilon q_1}{\tau}+q_2}{\sqrt{(1-\rho_2)}}+_o(\frac{1}{\tau}) \bigg) \\
  =&1-\Phi\bigg(\frac{-q_2}{\sqrt{(1-\rho_2)}}+\frac{\sqrt{\Phi_0}t+\mu_0-\sigma_\epsilon q_1}{\tau\sqrt{(1-\rho_2)}}+_o(\frac{1}{\tau}) \bigg)+ \\
  & \Phi\bigg(\frac{-q_2}{\sqrt{(1-\rho_2)}}-\frac{\sqrt{\Phi_0}t+\mu_0+\sigma_\epsilon q_1}{\tau\sqrt{(1-\rho_2)}}+_o(\frac{1}{\tau}) \bigg)
 \end{split}
 \end{equation}
 
 Therefore , type II error= 1-Power
 \begin{equation} \label{eqn10}
 \begin{split}
      & =\Phi\bigg(\frac{-q_2}{\sqrt{(1-\rho_2)}}+\frac{\sqrt{\Phi_0}t+\mu_0-\sigma_\epsilon q_1}{\tau\sqrt{(1-\rho_2)}}+_o(\frac{1}{\tau}) \bigg)- \\
  & \Phi\bigg(\frac{-q_2}{\sqrt{(1-\rho_2)}}-\frac{\sqrt{\Phi_0}t+\mu_0+\sigma_\epsilon q_1}{\tau\sqrt{(1-\rho_2)}}+_o(\frac{1}{\tau}) \bigg)
  \end{split}
 \end{equation}
 We now work with the test described by the  N.P. Lemma, which suggests that,
 Most powerful test (of its size) shall reject $H_{0i}$ if, ( suppose, $f_{H_{ji}}(y)$  is the pdf of Y under jth hypotheses ,j$\in \{0,1\}$ )
 \begin{equation*}
 \begin{split}
     & K(\tau)<  \frac{f_{H_{1i}}(y_i)}{f_{H_{0i}}(y_i)} \\
    \implies&K(\tau)< \frac{\frac{1}{\sqrt{2\pi\Phi_1(\tau)}}\exp\frac{-(y_i-\mu_1(\tau))^2}{2\Phi_1(\tau)}}{\frac{\tau}{\sqrt{2\pi\Phi_0}}\exp\frac{-(y_i-\frac{\mu_0}{\tau})^2\tau^2}{2\Phi_0}}\\
    \implies& K(\tau)< \exp\bigg[ -\frac{y_i^2}{2\Phi_1(\tau)}+\frac{y_i\mu_1(\tau)}{\Phi_1(\tau)}-\frac{\mu_1^2(\tau)}{2\Phi_1(\tau)}+\frac{y_i^2\tau^2}{2\Phi_0}-
     \frac{y_i\mu_0\tau}{\Phi_0}+\frac{\mu_0^2}{2\Phi_0}\\
     &\qquad-\log(\frac{\tau\sqrt{\Phi_1(\tau)}}{\sqrt{\Phi_0}})\bigg]\\
    \implies& \log(\frac{K(\tau)\tau\sqrt{\Phi_1(\tau)}}{\sqrt{\Phi_0}})< y_i^2\bigg(\frac{\tau^2}{2\Phi_0}-\frac{1}{2\Phi_1(\tau)}\bigg)-2y_i\bigg(\frac{\mu_0\tau}{2\Phi_0}-\frac{\mu_1(\tau)}{2\Phi_1(\tau)}\bigg)+\frac{\mu_0^2}{2\Phi_0}\\
  &\qquad-\frac{\mu_1^2(\tau)}{2\Phi_1(\tau)}
     \end{split}
 \end{equation*}
 From here, a simple algebra leads to the rejection region:

 \begin{multline}
     \bigg\{ y_i: y_i\bigg(\frac{\tau^2}{2\Phi_0}-\frac{1}{2\Phi_1(\tau)} \bigg)-\bigg(\frac{\mu_0\tau}{2\Phi_0}-\frac{\mu_1(\tau)}{2\Phi_1(\tau)}\bigg)\big>A(\tau)\bigg\}\\
     \mathlarger{\mathlarger{\mathlarger{\mathlarger{\cup}}}}\\
      \bigg\{ y_i: y_i\bigg(\frac{\tau^2}{2\Phi_0}-\frac{1}{2\Phi_1(\tau)} \bigg)-\bigg(\frac{\mu_0\tau}{2\Phi_0}-\frac{\mu_1(\tau)}{2\Phi_1(\tau)}\bigg)\big<-A(\tau)\bigg\}\\
     \end{multline}
      with
   \begin{multline*}
       A(\tau)=\\
       \bigg[\bigg(\frac{\tau^2}{2\Phi_0}-\frac{1}{2\Phi_1(\tau)}\bigg)\bigg(\log(\frac{K(\tau)\tau\sqrt{\Phi_1(\tau)}}{\sqrt{\Phi_0}})+\frac{\mu_1^2(\tau)}{2\Phi_1(\tau)}-\frac{\mu_0^2}{2\Phi_0}\bigg)+\bigg(\frac{\mu_0\tau}{2\Phi_0}-\frac{\mu_(\tau)}{2\Phi_1(\tau)}\bigg)^2\bigg]^{\frac{1}{2}}
   \end{multline*}
     i.e. rejection region is:
     \begin{equation}
         \{y_i:y_i>K_1(\tau)\}\cup\{y_i:y_i<K_2(\tau)\}
     \end{equation}
 
     where 
    
        \[
        K_1(\tau)= \frac{\bigg(\frac{\mu_0\tau}{2\Phi_0}-\frac{\mu_1(\tau)}{2\Phi_1(\tau)}\bigg)+A(\tau)}{\bigg(\frac{\tau^2}{2\Phi_0}-\frac{1}{2\Phi_1(\tau)} \bigg)}
        \]
    and
     \[
     K_2(\tau)= \frac{\bigg(\frac{\mu_0\tau}{2\Phi_0}-\frac{\mu_1(\tau)}{2\Phi_1(\tau)}\bigg)-A(\tau)}{\bigg(\frac{\tau^2}{2\Phi_0}-\frac{1}{2\Phi_1(\tau)} \bigg)}
     \]
   
   We shall get into the size and power calculations, but first concentrate on some approximations. This approximations are based on some previous calculations on $\mu_1(\tau), \Phi_1(\tau) \& \sqrt{\Phi_1(\tau)}.$
   
   \begin{equation*}
   \begin{split}
       \mu_1(\tau)=&q_2+\frac{\sigma_\epsilon q_1}{\tau}+_o(\frac{1}{\tau})\\
       \Phi_1(\tau)=&(1-\rho_2)+\frac{\Phi_0}{\tau^2}\\ \sqrt{\Phi_1(\tau)}=&\sqrt{(1-\rho_2)}+_o(\frac{1}{\tau}) 
   \end{split}
   \end{equation*}
   Therfore,
   \begin{equation*}
       \begin{split}
           (a)\frac{\mu_0\tau}{2\Phi_0}-\frac{\mu_1(\tau)}{2\Phi_1(\tau)}=&\frac{\mu_0\tau}{2\Phi_0}-\frac{q_2+\frac{\sigma_\epsilon q_1}{\tau}+_o(\frac{1}{\tau})}{2((1-\rho_2)+\frac{\Phi_0}{\tau^2})}\\
           =& \frac{\mu_0\tau}{2\Phi_0}-\frac{q_2}{2(1-\rho_2)}-\frac{\sigma_\epsilon q_1}{2(1-\rho_2)\tau}+_o(\frac{1}{\tau})\\
           (b)\frac{\tau^2}{2\Phi_0}-\frac{1}{2\Phi_1(\tau)}=&\frac{\tau^2}{2\Phi_0}-\frac{1}{2((1-\rho_2)+\frac{\Phi_0}{\tau^2})}\\
           =&\frac{\tau^2}{2\Phi_0}-\frac{1}{2(1-\rho_2)}+_o(\frac{1}{\tau})\\
           (c)\frac{\mu_1^2(\tau)}{2\Phi_1(\tau)}-\frac{\mu_0^2}{2\Phi_0}=&\frac{(q_2+\frac{\sigma_\epsilon q_1}{\tau}+_o(\frac{1}{\tau}))^2}{2((1-\rho_2)+\frac{\Phi_0}{\tau^2})}-\frac{\mu_0^2}{2\Phi_0}\\
           =&\frac{q_2^2+\frac{2\sigma_\epsilon q_1 q_2}{\tau}+_o(\frac{1}{\tau})}{2(1-\rho_2)+_o(\frac{1}{\tau})}-\frac{\mu_0^2}{2\Phi_0}\\
           =&\frac{q_2^2}{2(1-\rho_2)}-\frac{\mu_0^2}{2\Phi_0}+\frac{\sigma_\epsilon q_1 q_2}{\tau(1-\rho_2)}+_o(\frac{1}{\tau})
       \end{split}
   \end{equation*}
    Therefore,
   \begin{equation*}
       \begin{split}
           K_1(\tau)=&\frac{\bigg(\frac{\mu_0\tau}{2\Phi_0}-\frac{\mu_1(\tau)}{2\Phi_1(\tau)}\bigg)+A(\tau)}{\bigg(\frac{\tau^2}{2\Phi_0}-\frac{1}{2\Phi_1(\tau)} \bigg)}\\
           =&\frac{\bigg(\frac{\mu_0\tau}{2\Phi_0}-\frac{q_2}{2(1-\rho_2)}-\frac{\sigma_\epsilon q_1}{2(1-\rho_2)\tau}+_o(\frac{1}{\tau})\bigg)+A(\tau)}{\bigg(\frac{\tau^2}{2\Phi_0}-\frac{1}{2(1-\rho_2)}+_o(\frac{1}{\tau}) \bigg)}
       \end{split}
   \end{equation*}
   And,
   \begin{equation*}
       \begin{split}
           K_2(\tau)=&\frac{\bigg(\frac{\mu_0\tau}{2\Phi_0}-\frac{\mu_1(\tau)}{2\Phi_1(\tau)}\bigg)-A(\tau)}{\bigg(\frac{\tau^2}{2\Phi_0}-\frac{1}{2\Phi_1(\tau)} \bigg)}\\
           =&\frac{\bigg(\frac{\mu_0\tau}{2\Phi_0}-\frac{q_2}{2(1-\rho_2)}-\frac{\sigma_\epsilon q_1}{2(1-\rho_2)\tau}+_o(\frac{1}{\tau})\bigg)-A(\tau)}{\bigg(\frac{\tau^2}{2\Phi_0}-\frac{1}{2(1-\rho_2)}+_o(\frac{1}{\tau}) \bigg)}
       \end{split}
   \end{equation*}
  Where,
   \begin{multline*}
           A(\tau)=\\ \bigg[\bigg(\frac{\tau^2}{2\Phi_0}-\frac{1}{2\Phi_1(\tau)}\bigg)\bigg(\log(\frac{K(\tau)\tau\sqrt{\Phi_1(\tau)}}{\sqrt{\Phi_0}})+\frac{\mu_1^2(\tau)}{2\Phi_1(\tau)}-\frac{\mu_0^2}{2\Phi_0}\bigg)+\bigg(\frac{\mu_0\tau}{2\Phi_0}-\frac{\mu_(\tau)}{2\Phi_1(\tau)}\bigg)^2\bigg]^{\frac{1}{2}}\\
           =\bigg[\bigg(\frac{\tau^2}{2\Phi_0}-\frac{1}{2(1-\rho_2)}+_o(\frac{1}{\tau})\bigg)\bigg(\log(\frac{K(\tau)\tau\sqrt{\Phi_1(\tau)}}{\sqrt{\Phi_0}})+\frac{q_2^2}{2(1-\rho_2)}-\frac{\mu_0^2}{2\Phi_0}+ \\
       \frac{\sigma_\epsilon q_1 q_2}{\tau(1-\rho_2)}+_o(\frac{1}{\tau})\bigg)+\bigg(\frac{\mu_0\tau}{2\Phi_0}-\frac{q_2}{2(1-\rho_2)}-\frac{\sigma_\epsilon q_1}{2(1-\rho_2)\tau}+_o(\frac{1}{\tau})\bigg)^2\bigg]^{\frac{1}{2}}
   \end{multline*}
  
   Now, 
   \begin{equation}\label{eqn11}
       \begin{split}
          Size=& P_{H_{0i}}(\text{Rejection  Region})\\
          = & P_{H_{0i}}(Y>K_1(\tau))+P_{H_{0i}}(Y<K_2(\tau))\\
          = & P_{H_{0i}}(\frac{Y-\frac{\mu_0}{\tau}}{\frac{\sqrt{\Phi_0}}{\tau}}>\frac{K_1(\tau)-\frac{\mu_0}{\tau}}{\frac{\sqrt{\Phi_0}}{\tau}})+P_{H_{0i}}(\frac{Y-\frac{\mu_0}{\tau}}{\frac{\sqrt{\Phi_0}}{\tau}}<\frac{K_2(\tau)-\frac{\mu_0}{\tau}}{\frac{\sqrt{\Phi_0}}{\tau}})\\
          = & P(Z>\frac{\tau K_1(\tau)-\mu_0}{\sqrt{\Phi_0}})+P(Z<\frac{\tau K_2(\tau)-\mu_0}{\sqrt{\Phi_0}})\\
           & (\text{Where,} Z\sim N(0,1) \text{ and independent with Q}_1 \text{ and Q}_2 \text{ under H}_{0i})\\
          = & 1-\Phi\bigg(\frac{\tau K_1(\tau)-\mu_0}{\sqrt{\Phi_0}}\bigg)+\Phi\bigg(\frac{\tau K_2(\tau)-\mu_0}{\sqrt{\Phi_0}}\bigg)
       \end{split}
   \end{equation}
  Actual calculations need huge amount of labour and also do not simplify things at all. However, with the assumption $\tau\uparrow \infty$, a relatively simplified result can be obtained.\\ 
  For sufficiently large $\tau$,
  \begin{equation*}
      \begin{split}
          \tau K_1(\tau)\simeq &  \mu_0+\sqrt{2\Phi_0\big(\log(\frac{K(\tau)\tau\sqrt{\Phi_1(\tau)}}{\sqrt{\Phi_0}})+C_1\big)+\mu_0^2}\\
          \tau K_2(\tau)\simeq &  \mu_0-\sqrt{2\Phi_0\big(\log(\frac{K(\tau)\tau\sqrt{\Phi_1(\tau)}}{\sqrt{\Phi_0}})+C_1\big)+\mu_0^2}\\
      \end{split}
  \end{equation*}
  \[
  (C_1=\frac{q_2^2}{2(1-\rho_2)}-\frac{\mu_0^2}{2\Phi_0}+\frac{\sigma_\epsilon q_1 q_2}{\tau(1-\rho_2)})
  \]
  Therefore, equating size to $\alpha$, from \ref{eqn11} we obtain, for sufficiently large $\tau$,
  \[
  \alpha=1-\Phi(t_1)+\Phi(t_2)
  \]
  Where,
  \begin{equation*}
      \begin{split}
          t_1=&\frac{\tau K_1(\tau)-\mu_0}{\sqrt{\Phi_0}}\\
          \simeq & \frac{\sqrt{2\Phi_0\big(\log(\frac{K(\tau)\tau\sqrt{\Phi_1(\tau)}}{\sqrt{\Phi_0}})+C_1\big)+\mu_0^2}}{\sqrt{\Phi_0}}\\
          = &\sqrt{2\big(\log(\frac{K(\tau)\tau\sqrt{\Phi_1(\tau)}}{\sqrt{\Phi_0}})+C_1\big)+\frac{\mu_0^2}{\Phi_0}}
      \end{split}
  \end{equation*}
  and similarly,
   \begin{equation*}
      \begin{split}
          t_2=&\frac{\tau K_2(\tau)-\mu_0}{\sqrt{\phi_0}}\\
          \simeq &-\sqrt{2\big(\log(\frac{K(\tau)\tau\sqrt{\Phi_1(\tau)}}{\sqrt{\Phi_0}})+C_1\big)+\frac{\mu_0^2}{\Phi_0}}
      \end{split}
  \end{equation*}
  Resulting, $t_1=-t_2$, which implies, $t_1=z_{\frac{\alpha}{2}}. $ ($z_\gamma$ is the $100\gamma\%$ upper quantile of the standard normal distribution.)\\
   Therefore, for a $\alpha$ level of significance test, 
   \[
   \log(\frac{K(\tau)\tau\sqrt{\Phi_1(\tau)}}{\sqrt{\Phi_0}})=\frac{1}{2}\big[z_{\frac{\alpha}{2}}^2-\frac{\mu_0^2}{\Phi_0} \big]-C_1
   \]
  Accordingly the rejection region turns out to be:
  \[
  \{y:y> K_1^\alpha(\tau)\}\cup\{y:y< K_2^\alpha(\tau)\}
  \]
   with the expressions remaining same as before except
  \begin{equation*}
        \begin{split}
           &A^{\alpha}(\tau)\\
           =&\bigg[\bigg(\frac{\tau^2}{2\Phi_0}-\frac{1}{2(1-\rho_2)}+_o(\frac{1}{\tau})\bigg)\bigg(\frac{1}{2}\big[z_{\frac{\alpha}{2}}^2-\frac{\mu_0^2}{\Phi_0} \big]-C_1+\frac{q_2^2}{2(1-\rho_2)}-\frac{\mu_0^2}{2\Phi_0}+\\
       &\frac{\sigma_\epsilon q_1 q_2}{\tau(1-\rho_2)}+_o(\frac{1}{\tau})\bigg)+\bigg(\frac{\mu_0\tau}{2\Phi_0}-\frac{q_2}{2(1-\rho_2)}-\frac{\sigma_\epsilon q_1}{2(1-\rho_2)\tau}+_o(\frac{1}{\tau})\bigg)^2\bigg]^{\frac{1}{2}}\\
      =& \bigg[\bigg(\frac{\tau^2}{2\Phi_0}-\frac{1}{2(1-\rho_2)}+_o(\frac{1}{\tau})\bigg)\bigg(\frac{1}{2}\big[z_{\frac{\alpha}{2}}^2-\frac{\mu_0^2}{\Phi_0} \big]+_o(\frac{1}{\tau})\bigg)+\bigg(\frac{\mu_0\tau}{2\Phi_0}-\frac{q_2}{2(1-\rho_2)}-\\
      & \quad \frac{\sigma_\epsilon q_1}{2(1-\rho_2)\tau}+_o(\frac{1}{\tau})\bigg)^2\bigg]^{\frac{1}{2}}\\
 =&\bigg[\frac{\tau^2}{2\Phi_0}\bigg(\frac{1}{2}\big[z_{\frac{\alpha}{2}}^2-\frac{\mu_0^2}{\Phi_0} \big]\bigg)+\frac{\mu_0^2\tau^2}{4\Phi_0^2}+  _o(\tau^2)\bigg]^{\frac{1}{2}}\\
 =&\frac{\tau z_{\frac{\alpha}{2}}}{2\sqrt{\Phi_0}}+  _o(\tau)
 \end{split}
\end{equation*}

   Then,
   \begin{equation*}
       \begin{split}
           K_1^\alpha(\tau)=& \frac{\bigg(\frac{\mu_0\tau}{2\Phi_0}-\frac{q_2}{2(1-\rho_2)}-\frac{\sigma_\epsilon q_1}{2(1-\rho_2)\tau}+_o(\frac{1}{\tau})\bigg)+A^\alpha(\tau)}{\bigg(\frac{\tau^2}{2\Phi_0}-\frac{1}{2(1-\rho_2)}+_o(\frac{1}{\tau})\bigg)}\\
           =& \frac{\frac{(\mu_0+\sqrt{\Phi_0}z_{\frac{\alpha}{2}})}{2\Phi_0}\tau+  _o(\tau)}{\frac{\tau^2}{2\Phi_0}+_o(\tau^2)}\\
           =& \frac{(\mu_0+\sqrt{\Phi_0}z_{\frac{\alpha}{2}})}{\tau} +_o(\frac{1}{\tau})
       \end{split}
   \end{equation*}
     And similarly,  
      \begin{equation*}
       \begin{split}
           K_2^\alpha(\tau)=& \frac{\bigg(\frac{\mu_0\tau}{2\Phi_0}-\frac{q_2}{2(1-\rho_2)}-\frac{\sigma_\epsilon q_1}{2(1-\rho_2)\tau}+_o(\frac{1}{\tau})\bigg)-A^\alpha(\tau)}{\bigg(\frac{\tau^2}{2\Phi_0}-\frac{1}{2(1-\rho_2)}+_o(\frac{1}{\tau})\bigg)}\\
           =& \frac{\frac{(\mu_0-\sqrt{\Phi_0}z_{\frac{\alpha}{2}})}{2\Phi_0}\tau+  _o(\tau)}{\frac{\tau^2}{2\Phi_0}+_o(\tau^2)}\\
           =& \frac{(\mu_0-\sqrt{\Phi_0}z_{\frac{\alpha}{2}})}{\tau} +_o(\frac{1}{\tau})
       \end{split}
   \end{equation*}
   Therefore, test described by N.P.Lemma has rejection region of the form:
     
     \[
  \{\text{X}:\text{X}> K_1^{'\alpha}(\tau)\}\cup\{y:y< K_2^{'\alpha}(\tau)\}
  \]
     where 
    \begin{equation*}
       \begin{split}
           K_1^{'\alpha}(\tau)=& (\mu_0+\sqrt{\Phi_0}z_{\frac{\alpha}{2}}) +_o(1)
       \end{split}
   \end{equation*}
     And ,  
      \begin{equation*}
       \begin{split}
           K_2^{'\alpha}(\tau)=& (\mu_0-\sqrt{\Phi_0}z_{\frac{\alpha}{2}}) +_o(1)
       \end{split}
   \end{equation*}
     Which proves the second part of Result \ref{result2}.
     
   Therefore we get,
   \begin{equation} \label{eqn12}
       \begin{split}
           Power=& P_{H_{1i}}(Y> K_1^\alpha(\tau))+P_{H_{1i}}(Y<K_2^\alpha(\tau))\\
           =&P_{H_{1i}}\bigg(\frac{Y-\mu_1(\tau)}{\sqrt{\Phi_1(\tau)}}> \frac{K_1^\alpha(\tau)-\mu_1(\tau)}{\sqrt{\Phi_1(\tau)}}\bigg)+P_{H_{1i}}\bigg(\frac{Y-\mu_1(\tau)}{\sqrt{\Phi_1(\tau)}}< \frac{K_2^\alpha(\tau)-\mu_1(\tau)}{\sqrt{\Phi_1(\tau)}}\bigg)\\
           =&P\bigg(Z> \frac{K_1^\alpha(\tau)-\mu_1(\tau)}{\sqrt{\Phi_1(\tau)}}\bigg)+P\bigg(Z< \frac{K_2^\alpha(\tau)-\mu_1(\tau)}{\sqrt{\Phi_1(\tau)}}\bigg)\\
            & (\text{Where,} Z\sim N(0,1) \text{ and independent with Q}_1 \text{ and Q}_2 \text{ under H}_{1i})\\
           =& 1-\Phi\bigg( \frac{K_1^\alpha(\tau)-\mu_1(\tau)}{\sqrt{\Phi_1(\tau)}}\bigg)+\Phi\bigg( \frac{K_2^\alpha(\tau)-\mu_1(\tau)}{\sqrt{\Phi_1(\tau)}}\bigg)\\
           =& 1-\Phi\bigg( \frac{\frac{(\mu_0+\sqrt{\Phi_0}z_{\frac{\alpha}{2}}-\sigma_\epsilon q_1)}{\tau}-q_2}{\sqrt{(1-\rho_2)}) }+_o(\frac{1}{\tau})\bigg)+\Phi\bigg( \frac{\frac{(\mu_0-\sqrt{\Phi_0}z_{\frac{\alpha}{2}}-\sigma_\epsilon q_1)}{\tau}-q_2}{\sqrt{(1-\rho_2)} }+_o(\frac{1}{\tau})\bigg)\\
           =&1-\Phi\bigg(\frac{-q_2}{\sqrt{(1-\rho_2)}} +\frac{1}{\tau}\frac{(\mu_0+\sqrt{\Phi_0}z_{\frac{\alpha}{2}}-\sigma_\epsilon q_1)}{\sqrt{(1-\rho_2)}) }+_o(\frac{1}{\tau})\bigg)+\\
           & \Phi\bigg(\frac{-q_2}{\sqrt{(1-\rho_2)}} +\frac{1}{\tau}\frac{(\mu_0-\sqrt{\Phi_0}z_{\frac{\alpha}{2}}-\sigma_\epsilon q_1)}{\sqrt{(1-\rho_2)}) }+_o(\frac{1}{\tau})\bigg)
       \end{split}
   \end{equation}
  \ref{eqn12} yields type II error for the test based on N.P.Lemma =
   \begin{multline} \label{eqn13}
     TII_{NP}= \Phi\bigg(\frac{-q_2}{\sqrt{(1-\rho_2)}} +\frac{1}{\tau}\frac{(\mu_0+\sqrt{\Phi_0}z_{\frac{\alpha}{2}}-\sigma_\epsilon q_1)}{\sqrt{(1-\rho_2)}) }+_o(\frac{1}{\tau})\bigg)-\\
           \Phi\bigg(\frac{-q_2}{\sqrt{(1-\rho_2)}} +\frac{1}{\tau}\frac{(\mu_0-\sqrt{\Phi_0}z_{\frac{\alpha}{2}}-\sigma_\epsilon q_1)}{\sqrt{(1-\rho_2)}) }+_o(\frac{1}{\tau})\bigg)
   \end{multline}
   
   Recall, \ref{eqn10}:
    \begin{multline*}
    TII_{1}=\Phi\bigg(\frac{-q_2}{\sqrt{(1-\rho_2)}}+\frac{1}{\tau}\frac{(\mu_0+\sqrt{\Phi_0}t-\sigma_\epsilon q_1)}{\sqrt{(1-\rho_2)}}+_o(\frac{1}{\tau}) \bigg)- \\
\Phi\bigg(\frac{-q_2}{\sqrt{(1-\rho_2)}}-\frac{1}{\tau}\frac{(\mu_0+\sqrt{\Phi_0}t+\sigma_\epsilon q_1)}{\tau\sqrt{(1-\rho_2)}}+_o(\frac{1}{\tau}) \bigg)
 \end{multline*}
 where, t is the solution of:$\alpha=2-\Phi(t+\frac{2\mu_0}{\sqrt{\Phi_0}})-\Phi(t)$
 \paragraph{Note:}
 \begin{enumerate}
     \item The power we have calculated, are based on conditional test with respect to $Y = \frac{X_i}{\tau}|Q_1,Q_2$. To get a complete scenario, we need a measure based on the full data, which we obtain by taking expectation of the power (or equivalently , the type II error).
     \item  let, t be the solution of: $\alpha=2-\Phi(t+\frac{2\mu_0}{\sqrt{\Phi_0}})-\Phi(t)$. And let $t'$ be the solution of $\alpha=2-2\Phi(t')$. i.e. $t'=z_{\frac{\alpha}{2}}$. $\exists$ k ( because $\Phi(.)$ is a monotone function.) depending on ($q_1,q_2$) such that, t=$t'$+k. Now, if $\alpha\downarrow 0,$ $t\rightarrow \infty$ and $t'\rightarrow \infty$ with $\frac{t}{t'}\rightarrow 1$. Therefore, for sufficiently small $\alpha$, we can write, $t=t'+_o(t')$. i.e., $t=z_{\frac{\alpha}{2}}+_o(z_{\frac{\alpha}{2}})$

 \end{enumerate}
 We see, expanding through Taylor's series about $-\frac{q_2}{\sqrt{1-\rho_2}}$
 \begin{equation} \label{eqn14}
 \begin{split}
     TII_{NP}= & \Phi\bigg(\frac{-q_2}{\sqrt{(1-\rho_2)}} +\frac{1}{\tau}\frac{(\mu_0+\sqrt{\Phi_0}z_{\frac{\alpha}{2}}-\sigma_\epsilon q_1)}{\sqrt{(1-\rho_2)}) }+_o(\frac{1}{\tau})\bigg)-\\
          &  \Phi\bigg(\frac{-q_2}{\sqrt{(1-\rho_2)}} +\frac{1}{\tau}\frac{(\mu_0-\sqrt{\Phi_0}z_{\frac{\alpha}{2}}-\sigma_\epsilon q_1)}{\sqrt{(1-\rho_2)}) }+_o(\frac{1}{\tau})\bigg)\\
 = & \Phi\bigg(\frac{-q_2}{\sqrt{(1-\rho_2)}}\bigg)+\bigg(\frac{1}{\tau}\frac{\mu_0+\sqrt{\Phi_0}z_{\frac{\alpha}{2}}-\sigma_\epsilon q_1}{\sqrt{(1-\rho_2)}}+_o(\frac{1}{\tau}) \bigg)\phi\bigg(\frac{-q_2}{\sqrt{(1-\rho_2)}}\bigg)+_o(\frac{1}{\tau})-\\
 & \Phi\bigg(\frac{-q_2}{\sqrt{(1-\rho_2)}}\bigg)-\bigg(\frac{1}{\tau}\frac{\mu_0-\sqrt{\Phi_0}z_{\frac{\alpha}{2}}-\sigma_\epsilon q_1}{\sqrt{(1-\rho_2)}}+_o(\frac{1}{\tau}) \bigg)\phi\bigg(\frac{-q_2}{\sqrt{(1-\rho_2)}}\bigg)+_o(\frac{1}{\tau})\\
 =& \frac{1}{\tau}\frac{2\sqrt{\Phi_0}z_{\frac{\alpha}{2}}}{\sqrt{(1-\rho_2)}}\phi\bigg(\frac{-q_2}{\sqrt{(1-\rho_2)}}\bigg)+_o(\frac{1}{\tau})
 \end{split}
 \end{equation}
 Therefore,
 \begin{equation} \label{eqn15}
     \begin{split}
        & E_{Q_1,Q_2}( TII_{NP})\\
        = & \int_{-\infty}^\infty\int_{-\infty}^\infty \bigg(\frac{1}{\tau}\frac{2\sqrt{\Phi_0}z_{\frac{\alpha}{2}}}{\sqrt{(1-\rho_2)}}\phi\bigg(\frac{-q_2}{\sqrt{(1-\rho_2)}}\bigg)+_o(\frac{1}{\tau})\bigg)f(q_1)f(q_2)dq_1dq_2\\
         =&\frac{1}{\tau}\frac{2\sqrt{\Phi_0}z_{\frac{\alpha}{2}}}{\sqrt{(1-\rho_2)}} \int_{-\infty}^\infty \frac{1}{\sqrt{2\pi}}\exp\big[ -\frac{q_2^2}{2(1-\rho_2)}\big]\frac{1}{\sqrt{2\pi\rho_2}}\exp\big[ -\frac{q_2^2}{2\rho_2}\big]dq_2+_o(\frac{1}{\tau})\\
         =& \frac{1}{\tau}\frac{2\sqrt{\Phi_0}z_{\frac{\alpha}{2}}}{\sqrt{2\pi}}\int_{-\infty}^\infty \frac{1}{\sqrt{2\pi\rho_2(1-\rho_2)}}\exp\big[ -\frac{q_2^2}{2(1-\rho_2)\rho_2}\big]dq_2+_o(\frac{1}{\tau})\\
         =&\frac{1}{\tau}\frac{2\sqrt{\Phi_0}z_{\frac{\alpha}{2}}}{\sqrt{2\pi}}+_o(\frac{1}{\tau})
     \end{split}
 \end{equation}
 This proves 3rd part of the Result \ref{result2}\\
  Similarly, we get, again by  expanding through Taylor's series about $-\frac{q_2}{\sqrt{1-\rho_2}}$
  \begin{equation} \label{eqn16}
 \begin{split}
     TII_{1}= &\Phi\bigg(\frac{-q_2}{\sqrt{(1-\rho_2)}}+\frac{1}{\tau}\frac{\mu_0+\sqrt{\Phi_0}t-\sigma_\epsilon q_1}{\sqrt{(1-\rho_2)}}+_o(\frac{1}{\tau}) \bigg)- \\
 & \Phi\bigg(\frac{-q_2}{\sqrt{(1-\rho_2)}}-\frac{1}{\tau}\frac{\mu_0+\sqrt{\Phi_0}t+\sigma_\epsilon q_1}{\tau\sqrt{(1-\rho_2)}}+_o(\frac{1}{\tau}) \bigg)\\
 = & \Phi\bigg(\frac{-q_2}{\sqrt{(1-\rho_2)}}\bigg)+\bigg(\frac{1}{\tau}\frac{\mu_0+\sqrt{\Phi_0}t-\sigma_\epsilon q_1}{\sqrt{(1-\rho_2)}}+_o(\frac{1}{\tau}) \bigg)\phi\bigg(\frac{-q_2}{\sqrt{(1-\rho_2)}}\bigg)+_o(\frac{1}{\tau})-\\
 & \Phi\bigg(\frac{-q_2}{\sqrt{(1-\rho_2)}}\bigg)+\bigg(\frac{1}{\tau}\frac{\mu_0+\sqrt{\Phi_0}t-\sigma_\epsilon q_1}{\sqrt{(1-\rho_2)}}+_o(\frac{1}{\tau}) \bigg)\phi\bigg(\frac{-q_2}{\sqrt{(1-\rho_2)}}\bigg)+_o(\frac{1}{\tau})\\
 =& 2\bigg[\frac{1}{\tau}\frac{\sqrt{\Phi_0}t}{\sqrt{(1-\rho_2)}}+\frac{1}{\tau}\frac{\mu_0}{\sqrt{(1-\rho_2)}}\bigg]\phi\bigg(\frac{-q_2}{\sqrt{(1-\rho_2)}}\bigg)+_o(\frac{1}{\tau})
 \end{split}
 \end{equation}
 Now,
 \begin{equation}\label{eqn17}
     \begin{split}
         E_{Q_1,Q_2}(\mu_0\phi\bigg(\frac{-Q_2}{\sqrt{(1-\rho_2)}}\bigg) )=& \int_{-\infty}^\infty\int_{-\infty}^\infty (\sigma_\epsilon q_1+\sigma_0 q_2)\frac{1}{\sqrt{2\pi}}\exp\big[-\frac{q_2^2}{2(1-\rho_2)}\big]\times\\
         & \frac{1}{\sqrt{2\pi\rho_1}}\exp\big[-\frac{q_2^2}{2\rho_1}\big]\frac{1}{\sqrt{2\pi\rho_2}}\exp\big[-\frac{q_2^2}{2\rho_2}\big]dq_1dq_2
     \end{split}
 \end{equation}
 This can be expressed as sum of symmetric integrals of two odd functions. And hence both of them are 0, vanishing the whole term.i.e.
 \begin{equation*}
      E_{Q_1,Q_2}(\mu_0\phi\bigg(\frac{-Q_2}{\sqrt{(1-\rho_2)}}\bigg) )=0
 \end{equation*}
 Thus, \ref{eqn16}$\implies$
 \begin{equation} \label{eqn18}
     \begin{split}
         E(TII_{1})= & \int_{-\infty}^\infty\int_{-\infty}^\infty \bigg(\frac{1}{\tau}\frac{2\sqrt{\Phi_0}t}{\sqrt{(1-\rho_2)}}\phi\bigg(\frac{-q_2}{\sqrt{(1-\rho_2)}}\bigg)+_o(\frac{1}{\tau})\bigg)f(q_1)f(q_2)dq_1dq_2\\
         =& \frac{1}{\tau}\frac{2\sqrt{\Phi_0}}{\sqrt{(1-\rho_2)}}\int_{-\infty}^\infty\int_{-\infty}^\infty t\phi\bigg(\frac{-q_2}{\sqrt{(1-\rho_2)}}\bigg)f(q_1)f(q_2)dq_1dq_2+_o(\frac{1}{\tau})\\
         =& \frac{1}{\tau}\frac{2\sqrt{\Phi_0}}{\sqrt{(1-\rho_2)}}\int_{-\infty}^\infty\int_{-\infty}^\infty (z_{\frac{\alpha}{2}}+_o(z_{\frac{\alpha}{2}}))\phi\bigg(\frac{-q_2}{\sqrt{(1-\rho_2)}}\bigg)f(q_1)f(q_2)dq_1dq_2+_o(\frac{1}{\tau})\\
         =&\frac{1}{\tau}\frac{2\sqrt{\Phi_0}z_{\frac{\alpha}{2}}}{\sqrt{2\pi}}+_o(\frac{1}{\tau},z_{\frac{\alpha}{2}})
     \end{split}
 \end{equation}
 This proves the 4th part of Result \ref{result2}
   \end{proof}

\end{document}